\documentclass[reqno,11pt]{article}
\usepackage{mathrsfs}
\usepackage{amssymb}
\usepackage[usenames,dvipsnames]{xcolor}
\usepackage{amsthm}
\usepackage{amsmath}
\usepackage{amsfonts}
\usepackage{enumerate}
\usepackage{txfonts}
\usepackage{bm}

\usepackage{graphicx}

\numberwithin{equation}{section}

\newtheorem{theorem}{Theorem}[section]
\newtheorem{lemma}[theorem]{Lemma}
\newtheorem{corollary}[theorem]{Corollary}
\newtheorem{proposition}[theorem]{Proposition}

\theoremstyle{definition}
\newtheorem{definition}[theorem]{Definition}
\newtheorem{assumption}[theorem]{Assumption}
\newtheorem{example}[theorem]{Example}

\theoremstyle{remark}
\newtheorem{remark}[theorem]{Remark}

\begin{document}

\def\be{\begin{eqnarray}}
\def\ee{\end{eqnarray}}
\def\p{\partial}
\def\no{\nonumber}
\def\e{\epsilon}
\def\de{\delta}
\def\De{\Delta}
\def\om{\omega}
\def\Om{\Omega}
\def\f{\frac}
\def\th{\theta}
\def\la{\lambda}
\def\lab{\label}
\def\b{\bigg}
\def\var{\varphi}
\def\na{\nabla}
\def\ka{\kappa}
\def\al{\alpha}
\def\La{\Lambda}
\def\ga{\gamma}
\def\Ga{\Gamma}
\def\ti{\tilde}
\def\wti{\widetilde}
\def\wh{\widehat}
\def\ol{\overline}
\def\ul{\underline}
\def\Th{\Theta}
\def\si{\sigma}
\def\Si{\Sigma}
\def\oo{\infty}
\def\q{\quad}
\def\z{\zeta}
\def\co{\coloneqq}
\def\eqq{\eqqcolon}
\def\bt{\begin{theorem}}
\def\et{\end{theorem}}
\def\bc{\begin{corollary}}
\def\ec{\end{corollary}}
\def\bl{\begin{lemma}}
\def\el{\end{lemma}}
\def\bp{\begin{proposition}}
\def\ep{\end{proposition}}
\def\br{\begin{remark}}
\def\er{\end{remark}}
\def\bd{\begin{definition}}
\def\ed{\end{definition}}
\def\bpf{\begin{proof}}
\def\epf{\end{proof}}
\def\bex{\begin{example}}
\def\eex{\end{example}}
\def\bq{\begin{question}}
\def\eq{\end{question}}
\def\bas{\begin{assumption}}
\def\eas{\end{assumption}}
\def\ber{\begin{exercise}}
\def\eer{\end{exercise}}
\def\mb{\mathbb}
\def\mbR{\mb{R}}
\def\mbZ{\mb{Z}}
\def\mc{\mathcal}
\def\mcS{\mc{S}}
\def\ms{\mathscr}
\def\lan{\langle}
\def\ran{\rangle}
\def\lb{\llbracket}
\def\rb{\rrbracket}

\title{Liouville type theorems for the steady axially symmetric Navier-Stokes and magnetohydrodynamic equations}

\author{Dongho Chae\thanks{ Department of Mathematics, Chung-Ang University, Seoul 156-756, Republic of Korea. Email: dchae@cau.ac.kr.}\and {Shangkun Weng\thanks{ Pohang Mathematics Institute, Pohang University of Science and Technology. Pohang, Gyungbuk, 790-784, Republic of Korea. Email: skwengmath@gmail.com.}}}

\pagestyle{myheadings} \markboth{Liouville type theorems for steady axi-symmetric Navier-Stokes and MHD equations}{Liouville type theorems}\maketitle

\begin{abstract}
In this paper we study Liouville properties of smooth steady axially symmetric solutions of the Navier-Stokes equations.
First, we provide another version of the Liouville theorem of \cite{kpr15} in the case of zero swirl, where we replaced the Dirichlet integrability condition by mild decay conditions. Then we prove some Liouville theorems under the assumption $\|\f{u_r}{r}{\bf 1}_{\{u_r< -\f 1r\}}\|_{L^{3/2}(\mbR^3)}< C_{\sharp}$ where $C_{\sharp}$ is a universal constant to be specified. In particular, if $u_r(r,z)\geq -\f1r$ for $\forall (r,z)\in[0,\oo)\times\mbR$, then ${\bf u}\equiv 0$.  Liouville theorems also hold if $\displaystyle\lim_{|x|\to \oo}\Ga =0$ or $\Ga\in L^q(\mbR^3)$ for some $q\in [2,\oo)$ where $\Ga= r u_{\th}$. We also established some interesting inequalities for $\Om\co \f{\p_z u_r-\p_r u_z}{r}$, showing that $\na\Om$ can be bounded by $\Om$ itself. All these results are extended to the axially symmetric MHD and Hall-MHD equations with ${\bf u}=u_r(r,z){\bf e}_r +u_{\th}(r,z) {\bf e}_{\th} + u_z(r,z){\bf e}_z, {\bf h}=h_{\th}(r,z){\bf e}_{\th}$, indicating that the swirl component of the magnetic field does not affect the triviality. Especially, we establish the maximum principle for the total head pressure $\Phi=\f {1}{2} (|{\bf u}|^2+|{\bf h}|^2)+p$ for this special solution class.
\end{abstract}

\begin{center}
\begin{minipage}{5.5in}
\end{minipage}
\end{center}

\section{Introduction}\hspace*{\parindent}
The Steady Navier-Stokes takes the following form.
\be\lab{sns} &&\begin{cases}
({\bf u}\cdot\nabla) {\bf u} +\nabla p =\Delta {\bf u},\q \q \forall {\bf x}\in \mbR^3\\
\text{div }{\bf u}=0,
\end{cases}\\\lab{u-limit}
&&\displaystyle \lim_{|x|\to \oo} {\bf u}({\bf x}) =0.
\ee

Consider a weak solution to (\ref{sns}) with finite Dirichlet integral
\be\lab{ns-dirichlet}
\int_{\mbR^3} |\nabla {\bf u}({\bf x})|^2 d{\bf x}<+\oo.
\ee
It is well-known that a weak solution to (\ref{sns}) belonging to $W^{1,2}_{loc}(\mbR^3)$ is indeed smooth. Here we include the result stated in Galdi's book \cite{galdi11}.

\bt\lab{galdi}(Theorem X.5.1 in \cite{galdi11}).
{\it
Let ${\bf u}({\bf x})$ be a weak solution of (\ref{sns}) satisfying (\ref{ns-dirichlet}) and $p({\bf x})$ be the associated pressure, then there exists $p_1\in \mbR$ such that
\be\lab{asymptotic}
\lim_{|{\bf x}|\to \oo} |\na^{\al} {\bf u}({\bf x})| + \lim_{|{\bf x}|\to \oo} |\na^{\al} (p({\bf x})-p_1)|=0
\ee
uniformly for all multi-index $\al=(\al_1,\al_2,\al_3)\in [\mb{N}\cup \{0\}]^3$.
}
\et
We remark here that from (\ref{ns-dirichlet}) and (\ref{asymptotic}), one has ${\bf u}\in L^6(\mbR^3)\cap L^{\oo}(\mbR^3)$. So ${\bf u}\in L^p(\mbR^3)$ for any $p\in [6,\oo]$.

One of outstanding open problems for steady Navier-Stokes equations is: Is a smooth solution to (\ref{sns}) satisfying (\ref{ns-dirichlet}) identically zero? For the two dimensional case, the uniqueness result had been proved in \cite{gw78}. For the three dimensional case, Galdi \cite{galdi11} first established the triviality, under an additional integrability condition ${\bf u}\in L^{9/2}(\mbR^3)$. Chae and Yoneda \cite{cy13} also obtained a Liouville theorem if ${\bf u}\in X\cap Y$, where $X, Y$ are two function spaces and $X$ controls the high oscillation in large part of ${\bf u}$ and $Y$ gives control on the decay rate of ${\bf u}$ for sufficiently large ${\bf x}$. In \cite{chae14}, Chae showed that the condition $\Delta {\bf u}\in L^{6/5}(\mbR^3)$ is enough to guarantee the triviality. Note that this condition is stronger than the finite Dirichlet condition $\nabla {\bf u}\in L^2(\mbR^3)$, but both of them have the same scaling. In \cite{kpr15}, the authors proved that any axially symmetric smooth solution to (\ref{sns}) without swirl $u_{\th}\equiv 0$, but with finite Dirichelt integral must be zero.

To simplify the problem, in this article, we will consider the solution ${\bf u}$ with additional axially symmetric property. More precisely, we introduce the cylindrical coordinate
\be\no
r=\sqrt{x_1^2+x_2^2},\q \theta= \arctan\f{x_2}{x_1},\q z=x_3.
\ee
We denote ${\bf e}_r, {\bf e}_{\th}, {\bf e}_z$ the standard basis vectors in the cylindrical coordinate:
\be\no
{\bf e}_r= (\cos\th,\sin\th, 0),\q {\bf e}_{\th}= (-\sin\th,\cos\th,0),\q {\bf e}_z= (0,0,1).
\ee

A function $f$ is said to be {\it axially symmetric} if it does not depend on $\th$. A vector-valued function ${\bf u}= (u_r, u_{\th}, u_z)$ is called {\it axially symmetric} if $u_r, u_{\th}$ and $u_z$ do not depend on $\th$. A vector-valued function ${\bf u}= (u_r, u_{\th}, u_{z})$ is called {\it axially symmetric with no swirl} if $u_{\th}=0$ while $u_r$ and $u_z$ do not depend on $\th$. For more information about smooth axially symmetric vector fields, one may refer to \cite{lw09}.

Assume that ${\bf u}({\bf x})= u_r(r,z) {\bf e}_r+ u_{\th}(r,z) {\bf e}_{\th} + u_z(r,z){\bf e}_z$ is a smooth solution to (\ref{sns}). The corresponding asymmetric steady Navier-Stokes equations read as follows.
\be\lab{asym-sns}
\begin{cases}
(u_r\p_r+ u_z\p_z) u_r-\f{u_{\th}^2}{r} + \p_r p =\left(\p_r^2+\f{1}{r}\p_r+\p_z^2-\f{1}{r^2}\right) u_r, \\
(u_r\p_r+ u_z\p_z) u_{\th} +\f{u_r u_{\th}}{r}=\left(\p_r^2+\f{1}{r}\p_r+\p_z^2-\f{1}{r^2}\right)u_{\th},\\
(u_r\p_r+ u_z\p_z) u_z + \p_z p=\left(\p_r^2+\f{1}{r}\p_r+\p_z^2\right) u_z,\\
\p_r u_r +\f{u_r}{r} +\p_z u_z=0.
\end{cases}
\ee

In this paper, we want to investigate that what kinds of decay and integrability conditions we should prescribe on $u_r, u_{\th}$ and $u_{z}$ to guarantee the triviality of axially symmetric smooth solution to (\ref{sns}). Some of our conditions are weaker than the finite Dirichlet integral condition.

To our purpose, we also need the vorticity $\bm{\om}({\bf x})=\text{curl }{\bf u}({\bf x})= \om_r {\bf e}_r + \om_{\th} {\bf e}_{\th}+ \om_z {\bf e}_z$, where
\be\no
\om_r= -\p_z u_{\th},\q \om_{\th}= \p_z u_r -\p_r u_z,\q \om_z= \f 1r \p_r(r u_{\th}).
\ee
The equations satisfied by $\om_r,\om_{\th}$ and $\om_z$ are listed as follows.
\be\lab{vorticity}\begin{array}{ll}
(u_r\p_r+ u_z\p_z) \om_r - (\om_r\p_r +\om_z\p_z) u_r =\left(\p_r^2+\f{1}{r}\p_r+\p_z^2-\f{1}{r^2}\right)\om_r,\\
(u_r\p_r+ u_z\p_z) \om_{\th} - \f{u_r\om_{\th}}{r}-\f{1}{r}\p_z (u_{\th}^2) =\left(\p_r^2+\f{1}{r}\p_r+\p_z^2-\f{1}{r^2}\right)\om_{\th},\\
(u_r\p_r+ u_z\p_z) \om_z - (\om_r\p_r +\om_z\p_z) u_z =\left(\p_r^2+\f{1}{r}\p_r+\p_z^2\right)\om_z.
\end{array}\ee

This paper is structured as follows. In section \ref{lns}, we first provide another version of the Liouville theorem of \cite{kpr15} for axially symmetric flows with {\em no swirl}, but replacing the Dirichlet integral condition by mild decay condition at infinity. Then we establish some Liouville theorems under $\|\f{u_r}{r}{\bf 1}_{\{u_r< -\f 1r\}}\|_{L^{3/2}(\mbR^3)}\leq C_{\sharp}$, where $C_{\sharp}$ is a universal constant to be specified later. In particular, if $u_r(r,z)\geq -\f 1r$ for any $(r,z)\in [0,\oo)\times\mbR$, then ${\bf u}\equiv 0$. The Liouville theorem also holds under some decay or $L^q(\Bbb R^3)$ ($q\in [2,\oo)$) integrability conditions on $\Ga= ru_{\th}$. We also establish some interesting inequalities for $\Om$, use the equation for $\f{u_{\th}}{r}$ and the relation between $\f{u_r}{r}$ and $\Om$. These indicate that $\na\Om$ can be bounded by $\Om$ itself, which is not so clear at first sight due to the additional term $-\f{1}{r^2}\p_z(u_{\th})^2$. All these results will be extended to the axially symmetric MHD case in section \ref{lmhd}, where we consider a special solution form ${\bf u}({\bf x})= u_r(r,z){\bf e}_r+ u_{\th}(r,z){\bf e}_{\th}+ u_z(r,z){\bf e}_z$ and ${\bf h}({\bf x})= h_{\th}(r,z){\bf e}_{\th}$. Our results show that the swirl component of the magnetic field does not affect the triviality. Especially, we establish a maximum principle for the total head pressure $\Phi=\f{1}{2}(|{\bf u}|^2+ |{\bf h}|^2)+ p$. In section \ref{lhmhd}, we investigate the corresponding problem for steady resisitve, viscous Hall-MHD equations and obtain similar results as the MHD case.

\section{Liouville type theorems for steady Navier-Stokes equations}\lab{lns}

\subsection{The Liouville theorem for axially symmetric flows with no swirl}

In \cite{kpr15}, the authors have showed that the axially symmetric smooth solutions to (\ref{sns}) satisfying (\ref{ns-dirichlet}) must be zero in the absence of swirl. In the following, we provide another version of their Liouville theorem, replacing the Dirichlet integrability  condition by
mild decay conditions on ${\bf u}$. 

\bt\lab{liou-Omega}
{\it
Let ${\bf u}({\bf x})$ be an axially symmetric smooth solution to (\ref{sns})-(\ref{u-limit}) with no swirl. Then ${\bf u}\equiv 0$.
}
\et

Indeed this follows from Theorem 5.2 in \cite{knss} immediately. In Theorem 5.2 in \cite{knss}, they proved that for any bounded weak solution ${\bf u}$ of the unsteady Navier-Stokes equations in $\mbR^3\times (-\oo, 0)$, if ${\bf u}$ is axially symmetric with no swirl, then ${\bf u}= (0,0, b_3(t))$ for some bounded measurable function $b_3:(-\oo, 0)\to \mbR$. By (\ref{u-limit}), ${\bf u}\equiv 0$. Here we provide another proof of Theorem \ref{liou-Omega} by using the maximum principle of $\Om$.

\bpf[Proof of Theorem \ref{liou-Omega}]

Since ${\bf u}$ is smooth, by (\ref{u-limit}), one can conclude that $\na^{\al} {\bf u}\in L^{\oo}(\mbR^3)$ and
\be\lab{limit}
\displaystyle\lim_{|{\bf x}|\to \oo} |\na^{\al} {\bf u}({\bf x})| + \lim_{|{\bf x}|\to \oo} |\na^{\al} (p({\bf x})-p_1)|=0
\ee
uniformly for all multi-index $\al=(\al_1,\al_2,\al_3)\in [\mb{N}\cup \{0\}]^3$.

In the case $u_{\th}\equiv 0$, the equation for $\om_{\th}$ is reduced to
\be\lab{omegatheta}
(u_r\p_r+ u_z\p_z) \om_{\th} -\f{u_r\om_{\th}}{r} = \left(\p_r^2+\f{1}{r}\p_r + \p_z^2 -\f{1}{r^2}\right)\om_{\th}.
\ee
It is easy to derive the equation for $\Om$:
\be\lab{Omega}
(u_r\p_r+ u_z\p_z) \Om =\left(\p_r^2+\f 3r\p_r +\p_z^2\right)\Om.
\ee
As is well known, $(\p_r^2+\f 3r\p_r +\p_z^2)\Om$ can be written as the Laplacian $\De_5$ in $\mbR^5$ (see \cite{knss}). Write $r=\sqrt{y_1^2+y_2^2+y_3^2+y_4^2}$ and $y_5=z$, then (\ref{Omega}) becomes
\be\no
(u_r\p_r+ u_z\p_z) \Om =\Delta_5 \Om.
\ee
By (\ref{limit}), $\displaystyle\lim_{|{\bf x}|\to \oo} \Om({\bf x})=0$, one can employ the maximum principle arising from (\ref{Omega}) to show that $\Om\equiv 0$. Since $\text{curl }{\bf u}= \om_{\th}(r,z){\bf e}_{\th}$, we have $\text{curl }{\bf u}\equiv 0$.  Since  $\text{div }{\bf u}\equiv 0$, we have ${\bf u}({\bf x})=\nabla \phi({\bf x})$ for a harmonic function $\phi$.  The decay condition (\ref{u-limit}) implies ${\bf u}\equiv 0$ by the Liouville theorem for a harmonic function.

\epf

\subsection{Liouville type theorems conditioned on $u_r$ and $u_{\th}$.}\hspace*{\parindent}
We look at the equation $u_{\th}$ directly:
\be\lab{theta}
(u_r\p_r + u_z\p_z) u_{\th} + \f{u_r}{r} u_{\th} =(\p_r^2+ \f 1r\p_r +\p_z^2- \f{1}{r^2}) u_{\th}.
\ee

\bt\lab{u-r}
{\it Let ${\bf u}({\bf x})$ be an axially symmetric smooth solution to (\ref{sns})-(\ref{u-limit}). Suppose there exists $q\in[1, \infty)$ such that
\be\lab{utheta-q-2}
u_{\th}\in L^{q}(\mbR^3),
\ee
and
\be\lab{u-r1}
\|\f{u_r}{r}{\bf 1}_{\{u_r< -\f 1r\}}\|_{L^{3/2}(\mbR^3)}< \f{4(q+1)}{(q+2)^2 C_*^2},
\ee
where $C_*$ is the optimal constant in the Sobolev inequality $\|f\|_{L^6(\mbR^3)}\leq C_*\|\nabla f\|_{L^2(\mbR^3)}$ for any $f\in C_0^{\oo}(\mbR^3)$ and ${\bf 1}_{E}$ is the characteristic function for a set $E\subset \mbR^3$. Then, we have ${\bf u}\equiv 0$. In particular, if
\be\lab{u-r2}
u_r\geq  -\f{1}{r},\q \textit{for all $(r,z)\in [0,\oo)\times \mbR$},
\ee
then (\ref{u-r1}) is automatically satisfied and ${\bf u}\equiv 0$.
}
\et

\bpf
Since ${\bf u}$ is smooth, by (\ref{u-limit}) we have $\|{\bf u}\|_{L^{\oo}(\mbR^3)}<\oo$. 
Moreover, one can conclude that $\na^{\al} {\bf u}\in L^{\oo}(\mbR^3)$ and
\be\lab{up-limit}
\displaystyle\lim_{|{\bf x}|\to \oo} |\na^{\al} {\bf u}({\bf x})| + \lim_{|{\bf x}|\to \oo} |\na^{\al} (p({\bf x})-p_1)|=0
\ee
uniformly for all multi-index $\al=(\al_1,\al_2,\al_3)\in [\mb{N}\cup \{0\}]^3$. 

Introduce the radial cut-off function $\sigma\in C_0^{\oo}(\mathbb{R}^3)$ such that
\be\no
\si(x) =\si(|x|) =\begin{cases}
1\q\q \text{if } |x|<1,\\
0\q\q \text{if } |x|>2,
\end{cases}
\ee
and $0\leq \si(x)\leq 1$ for $1<|x|<2$. Without loss of generality, we may assume that $\si(|x|)$ is monotonic decreasing in $[0,+\oo)$. Then, for each $R>0$, we define $\si_R(x)= \si(|x|/R)$, then the support of $\nabla \si_R(x)$ is contained in $Q_R\co B_{2R}(0)\setminus B_R(0)$.

For $q\in [0,\oo)$, multiplying (\ref{theta}) by $ \si_R |u_{\th}|^{q} u_{\th}$, and integrating over $\mbR^3$, then we obtain
\be\lab{theta-q}
\int_{\mbR^3} \si_R |u_{\th}|^{q} u_{\th} (u_r\p_r+u_z\p_z+\f{u_r}{r}) u_{\th} dx = \int_{\mbR^3} \si_R |u_{\th}|^{q} u_{\th}\left(\p_r^2+\f 1r\p_r+\p_z^2-\f 1{r^2}\right) u_{\th} dx.
\ee
We estimate both sides as follows.
\be\no
LHS &=&\frac{1}{q+2}\int_{\mbR^3} \si_R(x) (u_r\p_r + u_z\p_z) |u_{\th}|^{q+2} dx + \int_{\mbR^3} \si_R(x) \f{u_r}{r} |u_{\th}|^{q+2} dx \\\no
&=& \f{2\pi}{q+2} \int_{-\oo}^{\oo} \int_0^{\oo}r \si\b(\f{\sqrt{r^2+z^2}}{R}\b)(u_r\p_r + u_z\p_z) |u_{\th}|^{q+2} dr dz +  \int_{\mbR^3} \si_R(x) \f{u_r}{r} |u_{\th}|^{q+2} dx\\\no
&=&-\f{2\pi}{q+2} \int_{-\oo}^{\oo} \int_0^{\oo}r \si'\b(\f{\sqrt{r^2+z^2}}{R}\b) |u_{\th}|^{q+2}\f{r u_r+z u_z}{R\sqrt{r^2+z^2}} drdz+ \int_{\mbR^3} \si_R(x) \f{u_r}{r} |u_{\th}|^{q+2} dx \\\no
&\co& A_1(R) +A_2(R).
\ee
\be\no
RHS &=& \int_{\mbR^3} \si_R(x) |u_{\th}|^{q} u_{\th} \left(\p_r^2+\f 1r\p_r +\p_z^2 -\f{1}{r^2}\right) u_{\th} dx \\\no
&=& 2\pi \int_{-\oo}^{\oo} \int_0^{\oo}r \si\b(\f{\sqrt{r^2+z^2}}{R}\b) |u_{\th}|^{q} u_{\th}\left(\p_r^2+\f{1}{r}\p_r + \p_z^2\right) u_{\th} dr dz - \int_{\mbR^3} \si_R(x) \f{|u_{\th}|^{q+2} }{r^2} dx \\\no
&=& -\f{8\pi (q+1)}{(q+2)^2} \int_{-\oo}^{\oo} \int_0^{\oo}r \si\b(\f{\sqrt{r^2+z^2}}{R}\b) (|\p_r |u_{\th}|^{\frac{q+2}{2}}|^2 + |\p_z |u_{\th}|^{\frac{q+2}{2}}|^2) dr dz\\\no
&\quad&+\f{2\pi}{q+2} \int_{-\oo}^{\oo} \int_0^{\oo}\f{1}{R}|u_{\th}|^{q+2} \b[\p_r\b(r\si'\b(\f{\sqrt{r^2+z^2}}{R}\b)\f{r}{\sqrt{r^2+z^2}}\b)+\p_z \b(r\si'\b(\f{\sqrt{r^2+z^2}}{R}\b)\f{z}{\sqrt{r^2+z^2}}\b)\b] dr dz\\\no
&\quad&- \int_{\mbR^3} \si_R(x) \f{|u_{\th}|^{q+2}}{r^2} dx \co B_1(R) +B_2(R) +B_3(R).
\ee

By the condition (\ref{u-limit}), we have (\ref{up-limit}) and
\be\no
\| {\bf u}\|_{L^\infty} + \|\f{u_{\th}}{r}\|_{L^{\oo}(\mbR^3)}<+\oo.
\ee
Since $\na |u_{\th}|^{\frac{q+2}{2}}= \f{q+2}{2} |u_{\th}|^{\frac{q}{2}-1}u_{\th}\na u_{\th}$ and $ u_{\th}\in L^\infty (\Bbb R^3)$,  (\ref{utheta-q-2}) implies $u_{\th}\in L^{q+2}(\mbR^3)$ and $\na |u_{\th}|^{\frac{q+2}{2}}\in L^2(\mbR^3)$. Now let $R$ tends to infinity, we obtain
\be\no
|A_1(R)|&\leq& \f{1}{q+2}(\|u_r\|_{L^{\oo}}+\|u_z\|_{L^{\oo}}) \f{1}{R}\int_{Q_R} |u_{\th}|^{q+2} dx \to 0,\q \text{as }R\to \oo,\\\no
A_2(R)&\to& \int_{\mbR^3} \f{u_r}{r} |u_{\th}|^{q+2} dx,\q \text{as } R\to \oo,\\\no
B_1(R) &=& -\f{4(q+1)}{(q+2)^2} \int_{\mbR^3} |\na |u_{\th}|^{\frac{q+2}{2}}|^2 dx,\q \text{as }R\to \oo,\\\no
|B_2(R)| &\leq& \f{C}{q+2}(\|\si'\|_{L^{\oo}}+\|\si^{''}\|_{L^{\oo}})\f{1}{R^2} \int_{Q_R} |u_{\th}|^{q+2} dx \to 0,\q \text{as }R\to \oo,\\\no
B_3(R) &=& -\int_{\mbR^3} \si_R(x)|\f{u_{\th}}{r}| |u_{\th}|^{q+2}dx\to -\int_{\mbR^3} \f{|u_{\th}|^{q+2}}{r^2} dx,\q \text{as }R\to \oo.
\ee

Combining all these calculations, we conclude that
\be\lab{q-integrability}
\f{4(q+1)}{(q+2)^2}\int_{\mbR^3} |\nabla |u_{\th}|^{\frac{q+2}{2}}|^2 dx +  \int_{\mbR^3} \f{|u_{\th}|^{q+2}}{r^2} dx + \int_{\mbR^3} \f{u_r}{r} |u_{\th}|^{q+2} dx =0.
\ee

If $\|\f{u_r}{r}{\bf 1}_{\{u_r< -\f 1r\}}\|_{L^{3/2}(\mbR^3)}<\f{4(q+1)}{(q+2)^2 C_*^2}$, then rewrite (\ref{q-integrability}) as
\be\no
&\quad&\f{4(q+1)}{(q+2)^2} \|\nabla |u_{\th}|^{\frac{q+2}{2}}\|_{L^2(\mbR^3)}^2+\int_{\mbR^3}\f{|u_{\th}|^{q+2}}{r^2} dx+ \int_{\mbR^3}\f{u_r}{r}{\bf 1}_{\{u_r\geq -\f 1r\}} |u_{\th}|^{q+2} dx \\\no
&=& -\int_{\mbR^3}\f{u_r}{r}|u_{\th}|^{q+2} {\bf 1}_{\{u_r< -\f 1r\}} dx\leq \|\f{u_r}{r}{\bf 1}_{\{u_r< -\f 1r\}}\|_{L^{3/2}(\mbR^3)}\||u_{\th}|^{\frac{q+2}{2}}\|_{L^6(\mbR^3)}^2 \\\no
&\leq& C_*^2\|\f{u_r}{r}{\bf 1}_{\{u_r< -\f 1r\}}\|_{L^{3/2}(\mbR^3)}\|\nabla |u_{\th}|^{\frac{q+2}{2}}\|_{L^2(\mbR^3)}^2<\f{4(q+1)}{(q+2)^2}\|\nabla |u_{\th}|^{\frac{q+2}{2}}\|_{L^2(\mbR^3)}^2,
\ee
yielding $\|\nabla |u_{\th}|^{\frac{q+2}{2}}\|_{L^2(\mbR^3)}^2=0$ and $u_{\th}\equiv 0$. Hence the problem reduces to the case of axially symmetric flows with no swirl, the theorem follows from Theorem \ref{liou-Omega}.

\epf

In the case of smooth axially symmetric solutions to (\ref{sns})-(\ref{u-limit}) with finite Dirichlet integral (\ref{ns-dirichlet}) one can have the optimal  constant in the estimate  (\ref{u-r1}).

\bt\lab{u-r-dirichlet}
{\it Let ${\bf u}({\bf x})$ be an axially symmetric smooth solution to (\ref{sns})-(\ref{u-limit}) with finite Dirichlet integral (\ref{ns-dirichlet}). Assume that
\be\lab{u-r3}
\|\f{u_r}{r}{\bf 1}_{\{u_r< -\f 1r\}}\|_{L^{3/2}(\mbR^3)}< \f{8}{9 C_*^2}.
\ee
Then ${\bf u}\equiv 0$. In particular, if
\be\lab{u-r4}
u_r\geq  -\f{1}{r},\q \textit{for all $(r,z)\in [0,\oo)\times \mbR$},
\ee
then (\ref{u-r3}) is automatically satisfied and ${\bf u}\equiv 0$.
}
\et

\bpf
Note that $f(q)\co \f{4(q+1)}{(q+2)^2}$ is monotonic decreasing for $q\in [1,\oo)$, and $\displaystyle\sup_{q\geq 1} \f{4(q+1)}{(q+2)^2}= f(1)=\f{8}{9}$. In the following, we will verify that (\ref{q-integrability}) holds for $q=1$. Since $\nabla {\bf u}\in L^2(\mbR^3)$, then by Sobolev embedding theorem ${\bf u}\in L^6(\mbR^3)$ and since
\be\no
\nabla {\bf u}&=& ({\bf e}_r \p_r + \f{{\bf e}_{\th}}{r}\p_{\th} + {\bf e}_z\p_z)(u_r {\bf e}_r + u_{\th} {\bf e}_{\th} + u_z {\bf e}_z)\\\no
&=& \p_r u_r {\bf e}_r \otimes {\bf e}_r + \f{u_r}{r} {\bf e}_{\th} \otimes {\bf e}_{\th} + \p_z u_r {\bf e}_z\otimes {\bf e}_r + \p_r u_{\th} {\bf e}_r\otimes {\bf e}_{\th} \\\no
&\q&-\f{u_{\th}}{r} {\bf e}_{\th}\otimes {\bf e}_r +\p_z u_{\th} {\bf e}_z\otimes {\bf e}_{\th} + \p_r u_z {\bf e}_r\otimes {\bf e}_z + \p_z u_z {\bf e}_z\otimes {\bf e}_z,
\ee
we also have $\|\nabla u_{\th}\|_{L^2(\mbR^3)}+\|\f{u_{r}}{r}\|_{L^2(\mbR^3)}+\|\f{u_{\th}}{r}\|_{L^2(\mbR^3)}\leq \|\nabla {\bf u}\|_{L^2(\mbR^3)}$. These implies that
\be\no
\b|\int_{\mbR^3} \f{u_r}{r} |u_{\th}|^3 dx\b| &\leq& \|\f{u_r}{r}\|_{L^2(\mbR^3)}\|u_{\th}\|_{L^6(\mbR^3)}<+\oo,\\\no
\|\nabla |u_{\th}|^{3/2}\|_{L^2(\mbR^3)}&\leq& C\|u_{\th}\|_{L^{\oo}(\mbR^3)}^{\f 12}\|\nabla {\bf u}\|_{L^2(\mbR^3)}<+\oo,\\\no
\b\|\f{|u_{\th}|^{\f 32}}{r}\b\|_{L^2(\mbR^3)} &\leq& \|u_{\th}\|_{L^{\oo}(\mbR^3)}^{\f 12} \|\f{u_{\th}}{r}\|_{L^2(\mbR^3)}<+\oo.
\ee

Letting $R\to \oo$ in (\ref{theta-q}), we obtain
\be\no
A_2&\to& \int_{\mbR^3} \f{u_r}{r} |u_{\th}|^{3} dx,\q\q \text{as } R\to \oo,  \\\no
B_1 &\to& -\f{8}{9}\int_{\mbR^3} |\nabla |u_{\th}|^{3/2}|^2 dx,\q\q \text{as } R\to \oo, \\\no
B_3 &\to& -\int_{\mbR^3} \f{|u_{\th}|^3}{r^2} dx,\q\q\text{as } R\to \oo,\\\no
|A_1|&\leq& C(\si) \|u_{\th}\|_{L^{6}(Q_R)}^3 \|(u_r, u_z)\|_{L^6(Q_R)},\q\text{as } R\to \oo,\\\no
|B_2|&\leq& C(\si) \|u_{\th}\|_{L^6(Q_R)}^2 R^{-1/2} \to 0,\q\text{as } R\to \oo.
\ee
Finally, we obtain (\ref{q-integrability}) for $q=1$.
\epf

\br\lab{u-r-remark}
{\it
It is well-known that if ${\bf u}$ is a smooth solution to (\ref{sns}), then so is ${\bf u}^{\la}({\bf x})= \la {\bf u}(\la {\bf x})$, for any $\la>0$. It should be emphasized here that the conditions (\ref{u-r3}),(\ref{u-r4})  and (\ref{u-r1}) are scaling invariant.
}
\er
Finally we give a simple vanishing criteria based on the $L^3(\mbR^3)$ integrability conditions for $u_r$ and $u_z$. Note that we do not need any additional conditions on $u_{\th}$.

\bt\lab{Only-u-rz}
{\it Let $({\bf u}, p)$ be an axially symmetric smooth solution to (\ref{sns})-(\ref{u-limit}) with finite integral (\ref{ns-dirichlet}). If $(u_r, u_z)\in L^3(\mbR^3)$, then ${\bf u}\equiv 0$.
}
\et

\bpf

Multiplying the first three equations in (\ref{asym-sns}) by $\si_R(x) u_r, \si_R(x) u_{\th}, \si_R(x) u_z$ and adding them together, and then integrating by parts, we obtain
\be\lab{Only-u-rz1}
&\q&\int_{\mbR^3} \si_R(x)\b(|\na u_r|^2+|\na u_{\th}|^2+ |\na u_z|^2+\f{u_r^2}{r^2}+\f{u_{\th}^2}{r^2}\b) dx\\\no
&=&\pi \int_{-\oo}^{\oo} \int_0^{\oo}\f{1}{R}|{\bf u}|^{2} \b[\p_r\b(r\si'\b(\f{\sqrt{r^2+z^2}}{R}\b)\f{r}{\sqrt{r^2+z^2}}\b)+\p_z \b(r\si'\b(\f{\sqrt{r^2+z^2}}{R}\b)\f{z}{\sqrt{r^2+z^2}}\b)\b] dr dz\\\no
&\q&+\int_{\mbR^3} (\f{1}{2}|{\bf u}|^2+ p) \f{1}{R}\si'\b(\f{\sqrt{r^2+z^2}}{R}\b)\f{r u_r+ z u_z}{\sqrt{r^2+z^2}} dx \co J_1 +J_2.
\ee
Since $\na {\bf u}\in L^2(\mbR^3)$, then $\f{1}{2}|{\bf u}|^2 +p \in L^3(\mbR^3)$ and the left hand side of (\ref{Only-u-rz1}) tends to
\be\no
\int_{\mbR^3}\b(|\na u_r|^2+|\na u_{\th}|^2+ |\na u_z|^2+\f{u_r^2}{r^2}+\f{u_{\th}^2}{r^2}\b) dx,\q \textit{as }R\to \oo.
\ee
If $(u_r,u_z)\in L^3(\mbR^3)$, then
\be\no
|J_1| &\leq& \f{C(\si)}{R}\||{\bf u}|\|_{L^3(Q_R)}^2 |Q_R|^{2/3}\to 0,\q\textit{as }R\to \oo,\\\no
|J_2| &\leq& \f{C(\si)}{R} \|\f{1}{2}|{\bf u}|^2+ p\|_{L^3(Q_R)}\|(u_r, u_z)\|_{L^3(Q_R)} |Q_R|^{1/3}\to 0,\q\textit{as }R\to \oo.
\ee
Hence letting $R\to \oo$ in (\ref{Only-u-rz1}), we can conclude that ${\bf u}\equiv 0$.

\epf

\subsection{Liouville type theorems conditioned on $\Ga= ru_{\th}$.}\hspace*{\parindent}
It follows from (\ref{asym-sns}) that
\be\lab{gamma-sns}
(u_r\p_r+ u_z\p_z)\Ga =\left(\p_r^2+\f{1}{r}\p_r+\p_z^2\right)\Ga -\f{2}{r}\p_r \Ga.
\ee
As is well known, for any harmonic function $h$ in $\mbR^3$, if $\lim_{|x|\to \oo} h(x)=0$ or $h\in L^q(\mbR^3)$ for $q\in [1,\oo)$, then $h\equiv 0$. Inspired by this, in the following we show that some decay or integrability conditions on $\Gamma$ are also enough to ensure that $\Gamma\equiv 0$.

\bt\lab{gamma condition}
{\it Let ${\bf u}({\bf x})$ be an axially symmetric smooth solution to (\ref{sns})-(\ref{u-limit}). If one of the following conditions holds
\begin{enumerate}[(i)]
  \item $\displaystyle \lim_{|{\bf x}|\to\oo} \Ga({\bf x})=0$,
  \item $\Ga\in L^{q}(\mbR^3)$ for some $q\in [2,\oo)$,
\end{enumerate}
then ${\bf u}\equiv 0$.
}
\et

\bpf  Assume (i), then the conclusion follows immediately from the maximum and minimum principle by $\Ga$. Note that $\Ga(0+,z)=0$ for $\forall z\in\mbR$, so $r=0$ does not cause any trouble. Now we assume (ii). Multiplying the equation (\ref{gamma-sns}) by $\si_R |\Ga|^{q-2} \Ga$ with $q\geq 2$, and then integrate over the whole space, we obtain that
\be\lab{integration-gamma-q}
\int_{\mbR^3} \si_R(x)|\Ga|^{q-2} \Ga (u_r\p_r+u_z\p_z) \Ga dx=  \int_{\mbR^3} \si_R(x) |\Ga|^{q-2} \Ga\b(\p_r^2+\f{1}{r}\p_r+\p_z^2-\f{2}{r}\p_r\b)\Ga dx.
\ee
We estimate both sides as follows.
\be\no
LHS&=&\f{2\pi}{q}\int_{-\oo}^{\oo}\int_0^{\oo} r\si\b(\f{\sqrt{r^2+z^2}}{R}\b) (u_r\p_r+u_z\p_z) |\Ga|^q(r,z) dr dz \\\no
&=&-\f{2\pi}{q}\int_{-\oo}^{\oo}\int_0^{\oo} r |\Ga|^q\si'\b(\f{\sqrt{r^2+z^2}}{R}\b)\f{r u_r+ z u_z}{R\sqrt{r^2+z^2}} dr dz\co -I(R),\\\no
RHS &=& 2\pi\int_{-\oo}^{\oo}\int_0^{\oo} r\si\b(\f{\sqrt{r^2+z^2}}{R}\b)|\Ga|^{q-2}\Ga(r,z)\b(\p_r^2+\f{1}{r}\p_r+\p_z^2-\f{2}{r}\p_r\b)\Ga(r,z) dr dz\\\no
&=& -\f{8\pi(q-1)}{q^2}\int_{-\oo}^{\oo}\int_0^{\oo} r\si\b(\f{\sqrt{r^2+z^2}}{R}\b)|\na |\Ga|^{q/2}|^2 dr dz+\f{4\pi}{q} \int_{-\oo}^{\oo}\int_0^{\oo} \si'\b(\f{\sqrt{r^2+z^2}}{R}\b)|\Ga|^{q} \f{r}{R\sqrt{r^2+z^2}}dr dz\\\no
&\quad&-\f{4\pi}{q} \int_{-\oo}^{\oo}\int_0^{\oo} r\si'\b(\f{\sqrt{r^2+z^2}}{R}\b)\f{r\p_r |\Ga|^{q}+z\p_z \p_r |\Ga|^{q}}{R\sqrt{r^2+z^2}} dr dz\\\no
&=& -\f{8\pi(q-1)}{q^2}\int_{-\oo}^{\oo}\int_0^{\oo} r\si\b(\f{\sqrt{r^2+z^2}}{R}\b)|\na |\Ga|^{q/2}|^2 dr dz+\f{4\pi}{q} \int_{-\oo}^{\oo}\int_0^{\oo} \si'\b(\f{\sqrt{r^2+z^2}}{R}\b)|\Ga|^{q} \f{r}{R\sqrt{r^2+z^2}}dr dz\\\no
&\quad&+\f{2\pi}{q}\int_{-\oo}^{\oo}\int_0^{\oo}\f{1}{R}|\Ga|^q \b[\p_r\b(r\si'\b(\f{\sqrt{r^2+z^2}}{R}\b)\f{r}{\sqrt{r^2+z^2}}\b)+\p_z \b(r\si'\b(\f{\sqrt{r^2+z^2}}{R}\b)\f{z}{\sqrt{r^2+z^2}}\b)\b] dr dz\\\no
&\co& -J_1(R) + J_2(R) +J_3(R).
\ee
Combining these together, we see that for any $0<R<R_1<\oo$:
\be\no
J_1(R) &=& I(R) + J_2(R) +J_3(R) \\\no
&\leq& J_1(R_1)= I(R_1) + J_2(R_1) + J_3(R_1)\\\lab{R}
&\leq& |I(R_1)| +|J_2(R_1)| + |J_3(R_1)|,
\ee
where we have used $J_1(R)$ is positive and monotonic increasing with respect to $R$.

Assume that $\Ga\in L^{q}(\mbR^3)$ for some $q\in [2,\oo)$, then by  H\"{o}lder's inequality,
\be\no
|I(R_1)|&\leq& \f{C(\si)}{R_1}\|(u_r,u_z)\|_{L^{\oo}(\mbR^3)} \|\Ga\|_{L^{q}(Q_{R_1})}\to 0,\q \q \text{as }R_1\to \oo,\\\no
|J_2(R_1)| &\leq& C R_1^{-2}\|\si'\|_{L^{\oo}}\int_{Q_{R_1}} |\Ga|^q  dx \leq C(\si) R_1^{-2}\|\Ga\|_{L^{q}(Q_{R_1})}\to 0,\q\q \text{as } R_1\to \oo,\\\no
|J_3(R_1)| &\leq& C R_1^{-2}(\|\si'\|_{L^{\oo}}+\|\si^{''}\|_{L^{\oo}})\int_{Q_{R_1}} |\Ga|^q  dx \leq C(\si) R_1^{-2}\|\Ga\|_{L^{q}(Q_{R_1})}\to 0,\q\q \text{as } R_1\to \oo.
\ee
In (\ref{R}), we first fix $R>0$ and let $R_1$ tends to $\oo$, we arrive at $\int_{B_R}|\nabla |\Ga|^{\f q2}|^2 dx=0$, which implies that $\nabla |\Ga|^{\f q2}\equiv 0$ in $B_R$. Since $R$ is arbitrary, we have $\nabla |\Ga|^{q/2}\equiv 0$ in the whole space. Since $\Ga\in L^q(\mbR^3)$ for some $q\in[2,\oo)$, then $\Ga\equiv 0$, i.e. $u_{\th}\equiv 0$.

The problem is reduced to the case of axially symmetric flows with no swirl, then the theorem follows from Theorem \ref{liou-Omega}. The proof is completed.

\epf

\br\lab{ga1}
{\it The condition in (i) is equivalent to $u_{\th}({\bf x})=o(\f 1r)$ as $|{\bf x}|\to \oo$. The bound $\f 1r$ is the optimal decay one can expect by the fundamental solution of the steady Stokes equations.}
\er
\br\lab{ga2}
{\it To understand the condition in (ii), we take $q=2$, then (ii) becomes $\Ga\in L^2(\mbR^3)$, i.e.
\be\lab{q2}
\int_0^{\oo} r^3 \int_{-\oo}^{\oo} |u_{\th}(r,z)|^2 dz dr<+\oo.
\ee
In \cite{cj09}, the authors have proved the following inequality by using $(\f{u_{\th}}{r}, \na u_{\th})\in L^2(\mbR^3)$:
\be\lab{cj}
\int_{-\oo}^{\oo} |u_{\th}(r,z)|^2 dz<+\oo,\q \forall r>0.
\ee
Comparing with (\ref{q2}) and (\ref{cj}), heuristically, one may guess that the integrability of $u_{\th}$ in the $z$-direction is enough, the decay rate of $u_{\th}$ in the radial direction maybe a key issue. Unfortunately, the decay rates obtained in \cite{cj09} seemed not good enough, this issue will be further investigated in \cite{weng151}.
}\er

From Theorem \ref{u-r}, it is natural to conject that a condition like $\|\f{u_{\th}}{r}\|_{X}\leq c_{\sharp}$, maybe enough to guarantee that $u\equiv 0$, where $X$ is a function norm and $c_{\sharp}>0$ is a constant depending only on the dimension $n=3$. Right now, we can not prove such a conjecture (see Corollary \ref{Omega-liouville}). However, we have some interesting conclusions.

\bt\lab{u-theta}
{\it
Let ${\bf u}({\bf x})$ be an axially symmetric smooth solution to (\ref{sns})-(\ref{u-limit}) with finite Dirichlet integral (\ref{ns-dirichlet}). If
\be\lab{theta condition}
\|\f{u_{\th}}{r}\|_{L^4(\mbR^3)}^2<\|\nabla \Om\|_{L^2(\mbR^3)},
\ee
then ${\bf u}\equiv 0$.
}
\et

\bpf
It follows from (\ref{vorticity}), then $\Om$ satisfies
\be\lab{generalOmega}
(u_r\p_r+ u_z\p_z)\Om -\f{1}{r^2} \p_z (u_{\th}^2)= \left(\p_r^2+\f{1}{r}\p_r+\p_z^2+\f{2}{r}\p_r\right)\Om.
\ee

We first verify that $\Om\in L^2(\mbR^3)$ and $\nabla\Om\in L^2(\mbR^3)$. Note that for any axially symmetric vector field ${\bf f}=(f_r, f_{\th}, f_z)$, we have
\be\no
|\na {\bf f}|^2&=& |\p_r f_r|^2+|\p_z f_r|^2+|\p_r f_{\th}|^2+|\p_z f_{\th}|^2+ \f{|f_r|^2}{r^2}\f{|f_{\th}|^2}{r^2}\\\no
&\q&+ |\p_r f_z|^2+|\p_z f_z|^2.
\ee
Hence by the definition of $\Om$, it suffices to verify that $\na^2 {\bf u}$ and $\na^3 {\bf u}$ belong to $L^2(\mbR^3)$. These indeed follow from a standard bootstrap argument by regarding $({\bf u}\cdot\nabla){\bf u}$ as a forcing term and using $L^p$ estimates for the Stokes system (see Theorem IV.2.1 in \cite{galdi11}). By using the cut-off function and integration by parts as in the proof of Theorem \ref{gamma condition}, we can conclude that
\be\lab{OmegaL2}
&\quad&\int_{\mbR^3} |\nabla\Om|^2 dx+ 2\pi \int_{-\oo}^{\oo} \Om^2(0,z) dz = \int_{\mbR^3}\f{1}{r^2}\p_z(u_{\th}^2) \Om dx\\\no
&=&- \int_{\mbR^3}\f{u_{\th}^2}{r^2} \p_z\Om dx\leq \|\f{u_{\th}}{r}\|_{L^4(\mbR^3)}^2 \|\nabla \Om\|_{L^2(\mbR^3)}\\\no
&<&  \|\nabla \Om\|_{L^2(\mbR^3)}^2,\q\q \text{if }(\ref{theta condition})\q \text{holds}.
\ee
Then $\Om\equiv 0$ if (\ref{theta condition}) holds. That is, $\om_{\th}= \p_z u_r-\p_r u_z\equiv 0$. Together with $\p_r u_r+\f{u_r}{r}+\p_z u_z=0$, we can conclude that $u_r=u_z\equiv 0$. The equation for $u_{\th}$ in (\ref{asym-sns}) reduces to
\be\no
\left(\p_r^2+\f{1}{r}\p_r+\p_z^2-\f{1}{r^2}\right) u_{\th}=0.
\ee
Setting $\La=\f{u_{\th}}{r}$, then
\be\no
\left(\p_r^2+\f{1}{r}\p_r+\p_z^2+\f{2}{r}\p_r\right) \La=0.
\ee
Same argument as in the proof of Theorem \ref{u-r}, we can show that $\La\equiv 0$ and $u_{\th}\equiv 0$.
\epf

\br\lab{zero-swirl}
{\it
This theorem shows that if the swirl component is smaller than the other two components in the sense of (\ref{theta condition}), then ${\bf u}\equiv 0$.
}
\er

In the following, we will use the equation for $\La\co \f{u_{\th}}{r}$ and the relation between $\f{u_r}{r}$ and $\Om$, to derive some interesting inequalities for $\Om$, showing that $\na\Om$ can be bounded by $\Om$ itself in some senses. Although $\Om$ satisfies an elliptic equation (\ref{generalOmega}), but due to the extra term $-\f{1}{r^2}\p_z (u_{\th}^2)$ in (\ref{generalOmega}), this property is not so clear at first sight.

\bt\lab{inverse}
{\it
Let ${\bf u}({\bf x})$ be an axially symmetric smooth solution to (\ref{sns})-(\ref{u-limit}) with finite Dirichlet integral (\ref{ns-dirichlet}). Then
\be\lab{inverse1}
\|\nabla \Om\|_{L^2(\mbR^3)}^2+2\pi\|\Om(0,\cdot)\|_{L^2(\mbR)}^2\leq 2^{5/3}C_* C_5^{2/3}\|\f{u_{\th}}{r}\|_{L^{\oo}} \|\f{u_{\th}}{r}\|_{L^2(\mbR^3)} \|\Om\|_{L^2(\mbR^3)}^{5/3}.
\ee
Furthermore, if we assume $\f{u_{\th}}{r}\in L^{3/2}(\mbR^3)$, then
\be\lab{inverse2}\begin{array}{ll}
&\|\nabla\Om\|_{L^2(\mbR^3)}^2+2\pi\|\Om(0,\cdot)\|_{L^2(\mbR)}^2\leq 2 C_*^4 C_6\|\f{u_{\th}}{r}\|_{L^{3/2}}^2\|\nabla\Om\|_{L^2}^2\\
&\quad\quad+ 2^{2/3} C_*^4\|\f{u_{\th}}{r}\|_{L^2}\|\f{u_{\th}}{r}\|_{L^{3/2}}\b(2\|\f{u_r}{r}\|_{L^3}+\|\p_z u_r\|_{L^3}+\|\p_z u_z\|_{L^3}\b)\|\Om\|_{L^2}^{2/3}\|\nabla\Om\|_{L^2}.
\end{array}\ee
}
\et

\bpf
It follows from (\ref{OmegaL2}) that
\be\lab{inverse11}
\|\nabla \Om\|_{L^2(\mbR^3)}^2 +2\pi \|\Om(0,\cdot)\|_{L^2(\mbR^3)}^2 &=&\b|-2\int_{\mbR^3}\La\p_z \La \cdot \Om dx \b|\\\no
&\leq& 2\|\La\|_{L^{\oo}}\|\p_z \La\|_{L^2(\mbR^3)}\|\Om\|_{L^2(\mbR^3)}.
\ee
It follows from (\ref{asym-sns}) that $\La$ satisfies the following equation
\be\lab{Lambda}
(u_r\p_r + u_z\p_z)\La +\f{2u_r}{r} \La= \left(\p_r^2+\f{1}{r}\p_r+\p_z^2+\f 2{r}\p_r\right)\La.
\ee
Since $\La\in L^2$ and $\nabla \La\in L^2(\mbR^3)$, then same as before, we can conclude
\be\lab{inverse12}
\|\nabla \La\|_{L^2(\mbR^3)}^2 + 2\pi\|\La(0,\cdot)\|_{L^2(\mbR)}^2 &=& -2\int_{\mbR^3}\f{u_r}{r} \La^2 dx\leq 2\|\f{u_r}{r}\|_{L^6(\mbR^3)}\|\La\|_{L^{12/5}(\mbR^3)}^2\\\no
&\leq& 2C_*^{3/2}\b\|\nabla\b(\f{u_r}{r}\b)\b\|_{L^2(\mbR^3)}\|\La\|_{L^2(\mbR^3)}^{3/2}\|\nabla\La\|_{L^2(\mbR^3)}^{1/2},
\ee
where we have used the classical Gagliardo-Nirenberg inequality on $\mbR^n$:
\be\no
\|f\|_{L^q}\leq C\|\nabla f\|_{L^p}^{\alpha} \|f\|_{L^s}^{1-\alpha},
\ee
with $p, q, s\geq 1$ and $\al\in [0,1]$ satisfy the identity
\be\no
\alpha=\b(\f1s-\f1q\b)\b(\f1n-\f1p+\f1s\b)^{-1}.
\ee

Then we obtain
\be\lab{inverse13}
\|\nabla \La\|_{L^2(\mbR^3)}\leq 2^{2/3}C_* \b\|\nabla\b(\f{u_r}{r}\b)\b\|_{L^2(\mbR^3)}^{2/3}\|\La\|_{L^2(\mbR^3)}.
\ee
We still need to explore the relation between $\f{u_r}{r}$ and $\Om$ by introducing the stream function, this relation was already known in the unsteady case \cite{cweng,lei}. By the divergence free condition, $\p_r (r u_r)+\p_z (r u_z)=0$, one can introduce a stream function $\psi_{\theta}$ such that
\be\no
u_r =-\p_z \psi_{\theta},\quad u_z=\f{1}{r}\p_r(r\psi_{\theta}).
\ee
Since $\omega_{\theta}=\p_z u_r-\p_r u_z$, we have
\be\no
-(\p_r^2 +\f{1}{r}\p_r+\p_z^2-\f{1}{r^2})\psi_{\theta}= \omega_{\theta}.
\ee
Setting $\varphi=\f{\psi_{\theta}}{r}$, then it is easy to see that
\be\no
-(\p_r^2+\f{3}{r}\p_r+\p_z^2)\varphi=\Omega.
\ee
The second order operator $(\p_r^2+\f{3}{r}+\p_z^2)$ can be interpreted as the Laplace operator in $\mbR^5$, see \cite{hll,knss}. Introduce
\be\no
y=(y_1,y_2,y_3,y_4, z),\quad r=\sqrt{y_1^2+y_2^2+y_3^2+y_4^2},\quad \Delta_y=(\p_r^2+\f{3}{r}\p_r+\p_z^2).
\ee
Hence we have $\varphi=(-\Delta_y)^{-1} \Omega$ and $\f{u_r}{r}=-\p_z \varphi$. By simple calculations, one has
\be\no
|\nabla_y^2 \varphi|^2\simeq |\p_r^2 \varphi|^2+|\f{1}{r}\p_r \varphi|^2+|\p_z^2 \varphi|^2+|\p_{rz}^2 \varphi|^2
\ee
and
\be\no
\int_{\mbR^3} |\nabla^2 \varphi|^2 dx&\leq& C_3\int_{-\infty}^{\infty} \int_0^{\infty} \b(|\p_r^2 \varphi|^2+|\f{1}{r}\p_r \varphi|^2+|\p_z^2 \varphi|^2+|\p_{rz}^2 \varphi|^2\b) r dr dz\\\no
&=&C_3\int_{-\infty}^{\infty} \int_0^{\infty} \b(|\p_r^2 \varphi|^2+|\f{1}{r}\p_r \varphi|^2+|\p_z^2 \varphi|^2+|\p_{rz}^2 \varphi|^2\b)w(r) r^3 dr dz\\\no
&\leq& C_4\int_{-\infty}^{\infty} \int_0^{\infty} |\nabla_y^2 \varphi|^2w(r) r^3 dr dz=C_4 \int_{-\infty}^{\infty} \int_0^{\infty} |\nabla_y^2 (-\Delta_y)^{-1}\Omega|^2w(r) r^3 dr dz\\\no
&=& C_4\int_{\mbR^5}|\nabla_y^2 (-\Delta_y)^{-1}\Omega|^2 w(r) dy\\\lab{inverse14}
&\leq& C_5\int_{\mbR^5} |\Omega|^2 w(r) dy= C_5\int_{\mbR^3} |\Omega|^2 dx,
\ee
where $w(r)=r^{-2}$ and in the last step we have used the boundedness of Riesz operators in weighted Sobolev spaces (Lemma 2 in \cite{hll}). See also Corollary 2 in \cite{cl02} for a similar weighted estimate for a singular integral operator.

Similarly, we also have
\be\lab{inverse15}
\int_{\mbR^3} |\nabla^2 \p_z\varphi|^2 dx\leq C_6\int_{\mbR^3} |\p_z \Omega|^2 dx.
\ee
Then by (\ref{inverse14}), we have
\be\lab{inverse16}
\b\|\nabla\b(\f{u_r}{r}\b)\b\|_{L^2(\mbR^3)}\leq C_5\|\Om\|_{L^2(\mbR^3)},\q \b\|\nabla\p_z\b(\f{u_r}{r}\b)\b\|_{L^2(\mbR^3)}\leq C_6\|\p_z\Om\|_{L^2(\mbR^3)}.
\ee
Then (\ref{inverse1}) follows from (\ref{inverse11}), (\ref{inverse13}) and (\ref{inverse16}).

To derive (\ref{inverse2}), we estimate (\ref{inverse11}) as follows.
\be\no
\|\nabla\Om\|_{L^2(\mbR^3)}^2+ 2\pi\|\Om(0,\cdot)\|_{L^2(\mbR)}^2 &\leq&\|\La\|_{L^{3/2}(\mbR^3)}\|\p_z\La\|_{L^6(\mbR^3)}\|\Om\|_{L^6(\mbR^3)}\\\lab{inverse21}
&\leq& C_*^2\|\La\|_{L^{3/2}(\mbR^3)}\|\nabla\p_z\La\|_{L^2(\mbR^3)}\|\nabla \Om\|_{L^2(\mbR^3)}.
\ee
To estimate $\nabla \p_z\La$, we first derive the equation for $\p_z\La$:
\be\lab{inverse22}
\left(\p_r^2+\f1r\p_r+ \p_z^2+\f{2}{r}\p_r\right)\p_z\La=(u_r\p_r+u_z\p_z)\p_z\La+2\f{u_r}{r}\p_z\La+2\p_z\b(\f{u_r}{r}\b)\La+ (\p_z u_r\p_r \La+\p_z u_z\p_z\La).
\ee
Then by integrate by parts, we have
\be\no
&\quad&\|\nabla\p_z\La\|_{L^2}^2+2\pi\|\p_z\La(0,\cdot)\|_{L^2(\mbR)}^2\\\no
&=&\b|\int_{\mbR^3}\b[2\f{u_r}{r}(\p_z \La)^2+ 2\p_z\b(\f{u_r}{r}\b)\La\p_z\La+ \p_z u_r\p_r \La\p_z\La +\p_z u_z(\p_z\La)^2\b]dx\b|\\\no
&\leq& 2\|\f{u_r}{r}\|_{L^3}\|\p_z\La\|_{L^2}\|\p_z\La\|_{L^6}+ 2\b\|\p_z\b(\f{u_r}{r}\b)\b\|_{L^6}\|\La\|_{L^{3/2}}\|\p_z\La\|_{L^6}\\\no
&\quad&\quad+ \|\p_z u_r\|_{L^3}\|\p_r\La\|_{L^2}\|\p_z\La\|_{L^6}+ \|\p_z u_z\|_{L^3}\|\p_z\La\|_{L^2}\|\p_z\La\|_{L^6}\\\no
&\leq& 2C_*\|\f{u_r}{r}\|_{L^3}\|\p_z\La\|_{L^2}\|\nabla\p_z\La\|_{L^2}+ 2C_*^2\b\|\nabla\p_z\b(\f{u_r}{r}\b)\b\|_{L^2}\|\La\|_{L^{3/2}}\|\nabla\p_z\La\|_{L^2}\\\lab{inverse23}
&\quad&\quad+ C_*\|\p_z u_r\|_{L^3}\|\p_r\La\|_{L^2}\|\nabla\p_z\La\|_{L^2}+ C_*\|\p_z u_z\|_{L^3}\|\p_z\La\|_{L^2}\|\nabla\p_z\La\|_{L^2},
\ee
By (\ref{inverse15}), (\ref{inverse13}), (\ref{inverse23}) yields
\be\lab{inverse24}\begin{array}{ll}
\|\nabla\p_z\La\|_{L^2}&\leq 2C_*^2C_6\|\La\|_{L^{3/2}}\|\nabla\Om\|_{L^2}+ C_*\b(2\|\f{u_r}{r}\|_{L^3}+\|\p_z u_r\|_{L^3}+\|\p_z u_z\|_{L^3}\b)\|\nabla\La\|_{L^2}\\
&\leq2C_*^2C_6\|\La\|_{L^{3/2}}\|\nabla\Om\|_{L^2}+ 2^{2/3} C_*^2\|\La\|_{L^2}\b(2\|\f{u_r}{r}\|_{L^3}+\|\p_z u_r\|_{L^3}+\|\p_z u_z\|_{L^3}\b)\|\Om\|_{L^2}^{2/3}.
\end{array}\ee
Hence (\ref{inverse2}) follows from (\ref{inverse21}) and (\ref{inverse24}).
\epf

\bc\lab{Omega-liouville}
{\it
Let ${\bf u}({\bf x})$ be an axially-symmetric smooth solution to (\ref{sns})-(\ref{u-limit}) with finite Dirichlet integral (\ref{ns-dirichlet}). If
\be\lab{u-theta}\begin{cases}
\|\f{u_{\th}}{r}\|_{L^{3/2}}^2<\f 1{4 C_*^4 C_6}\\
2^{2/3} C_*^4\|\f{u_{\th}}{r}\|_{L^2}\|\f{u_{\th}}{r}\|_{L^{3/2}}\b(2\|\f{u_r}{r}\|_{L^3}+\|\p_z u_r\|_{L^3}+\|\p_z u_z\|_{L^3}\b)\|\Om\|_{L^2}^{2/3}<\f 12\|\nabla\Om\|_{L^2},
\end{cases}\ee
then $\Om\equiv 0$, which implies that ${\bf u}\equiv 0$.
}
\ec

\br\lab{Omega-remark}
{\it
We remark here all the quantities in (\ref{inverse1}) and (\ref{inverse2}) have same scaling properties.
}\er


\section{Liouville type theorem for steady MHD equations}\lab{lmhd}\hspace*{\parindent}
In this section, we will investigate the steady Magnetohydrodynamics equations (MHD). The steady MHD equations are listed as follows.
\be\label{smhd}\begin{cases}
({\bf u}\cdot\nabla) {\bf u}+\nabla p =({\bf h}\cdot\nabla) {\bf h}+ \Delta {\bf u},\q\q\forall x\in \mbR^3,\\
\text{div } {\bf u}=0,\\
({\bf u}\cdot\nabla){\bf h}-({\bf h}\cdot\nabla){\bf u}=\Delta {\bf h},\q\q\forall x\in \mbR^3,\\
\text{div } {\bf h}=0,\\
\displaystyle\lim_{|x|\to \oo} {\bf u}({\bf x}) =\lim_{|x|\to \oo} {\bf h}({\bf x})= {\bf 0}.
\end{cases}
\ee
We also consider the weak solution to (\ref{smhd}) with finite Dirichlet integral:
\be\lab{mhd-dirichlet}
\int_{\mbR^3} |\na {\bf u}({\bf x})|^2 +|\na {\bf h}({\bf x})|^2 d{\bf x} <\oo.
\ee
Then following the argument developed in \cite{galdi11}, one can show that any weak solution to (\ref{smhd}) satisfying (\ref{mhd-dirichlet}) is smooth. A natural problem is whether this solution is zero or not. We first present the following simple Liouville theorem for MHD equations, by prescribing the $L^3$ integrability only on ${\bf u}$. Similar regularity criteria has been developed in \cite{hx1} and \cite{zhou} for unsteady MHD equations.
\bt\lab{Liou-mhd}
{\it Let $({\bf u}, {\bf h})$ be a smooth solution to (\ref{smhd}) in $\mbR^3$ with finite Dirichlet integral. If ${\bf u}\in L^3(\mbR^3)$, then ${\bf u}= {\bf h}\equiv 0$.
}
\et

\bpf
Multiplying the first and third equation in (\ref{smhd}) by $\si_R {\bf u}$ and $\si_R {\bf h}$ respectively, and integrating by parts, we finally obtain
\be\lab{integration-mhd}\begin{array}{ll}
&\quad\int_{\mbR^3} \si_R(x) (|\nabla {\bf u}|^2 +|\nabla {\bf h}|^2) dx + \f{1}{R} \int_{\mbR^3} (u_i\p_j u_i+ h_i\p_j h_i)\p_j \si(x/R) dx\\
&-\int_{\mbR^3}\f{1}{R}(\f12 (|{\bf u}|^2+|{\bf h}|^2)+ p) {\bf u}\cdot\nabla \si(x/R) dx + \int_{\mbR^3}\f{1}{R} ({\bf u}\cdot {\bf h}) ({\bf h}\cdot \nabla) \si(x/R) dx \\
&\co H_1 +H_2 +H_3 +H_4=0.
\end{array}\ee
It is well-known that ${\bf u}\in L^6(\mbR^3), {\bf h}\in L^6(\mbR^3)$ and $p\in L^3(\mbR^3)$. Since we assume that ${\bf u}\in L^q(\mbR^3)$ for $q\in [1,3]$, then ${\bf u}\in L^3(\mbR^3)$. It is easy to see that $H_1\to \|\nabla {\bf u}\|_{L^2}^2+ \|\nabla {\bf h}\|_{L^2}^2$ as $R\to \oo$, and also
\be\no
|H_2|&\leq& \|{\bf u}\|_{L^6(Q_R)} \|\nabla {\bf u}\|_{L^2(Q_R)}\to 0,\q \text{as } R\to \oo,\\\no
|H_3|&\leq& (\|{\bf u}\|_{L^6(Q_R)}^2+\|{\bf h}\|_{L^6(Q_R)}^2)\|{\bf u}\|_{L^3(Q_R)} \to 0,\q \text{as }R\to \oo,\\\no
|H_4|&\leq& \|{\bf u}\|_{L^3(Q_R)}\|{\bf h}\|_{L^6(Q_R)}^2 \to 0,\q \text{as }R\to \oo.
\ee
Hence letting $R\to \oo$ in (\ref{integration-mhd}), we obtain
\be\no
\|\nabla {\bf u}\|_{L^2(\mbR^3)}^2 + \|\nabla {\bf h}\|_{L^2(\mbR^3)}^2=0
\ee
Hence ${\bf u}={\bf h}\equiv 0$.
\epf

In the following, we derive the equation for the total head pressure $\Phi\co \f 12 (|{\bf u}|^2+|{\bf h}|^2)+ p$. Assume that $({\bf u}, {\bf h})$ is a smooth solution to steady MHD equations (\ref{smhd}), then by simple calculations, we obtain
\be\lab{pressure eq}
\Delta p &=& -\text{div }[({\bf u}\cdot\nabla){\bf u}]+ \text{div }[({\bf h}\cdot\nabla){\bf h}] \\\no
&=&- \sum_{i,j=1}^3 \p_i u_j \p_j u_i + \sum_{i,j=1}^3 \p_i h_j \p_j h_i.
\ee
This yields
\be\lab{bernoulli eq}
&\quad&({\bf u}\cdot\nabla)\Phi -({\bf h}\cdot \nabla)({\bf u}\cdot {\bf h})\\\no
&=& {\bf u}\cdot \Delta {\bf u}+ {\bf h}\cdot \Delta {\bf h}= \Delta\b(\f 12 (|{\bf u}|^2+|{\bf h}|^2)\b)- |\nabla {\bf u}|^2- |\nabla {\bf h}|^2 \\\no
&=& \Delta \Phi- \b(|\nabla {\bf u}|^2-\sum_{i,j=1}^3 \p_i u_j\p_j u_i\b)-\b(|\nabla {\bf h}|^2+ \sum_{i,j=1}^3 \p_i h_j \p_j h_i\b).
\ee
Since
\be\no
|\nabla {\bf u}|^2-\sum_{i,j=1}^3 \p_i u_j\p_j u_i&=& (\p_2 u_3 -\p_3 u_2)^2 +(\p_3 u_1-\p_1 u_3)^2 + (\p_1 u_2-\p_2 u_1)^2= |\text{curl }{\bf u}|^2\\\no
|\nabla {\bf h}|^2+\sum_{i,j=1}^3 \p_i h_j\p_j h_i&=&2\sum_{i=1}^3(\p_i h_i)^2 + (\p_2 h_3+\p_3 h_2)^2+ (\p_3 h_1+\p_1 h_3)^2 + (\p_1 h_2+\p_2 h_1)^2,
\ee
we have
\be\lab{bernoulli eq1}
({\bf u}\cdot\nabla) \Phi- ({\bf h}\cdot\nabla)({\bf u}\cdot {\bf h})-\Delta \Phi \leq 0.
\ee

Consider the special case, where $({\bf u}, {\bf h})$ are axi-symmetric, and of the following special form
\be\lab{asymmetric}
\begin{cases}
{\bf u}({\bf x})= u_r(r,z){\bf e}_r + u_{\th}(r,z) {\bf e}_{\th}+ u_z(r,z){\bf e}_z,\\
{\bf h}({\bf x})= h_{\th}(r,z) {\bf e}_{\th},
\end{cases}
\ee
then $({\bf h}\cdot \nabla) ({\bf u}\cdot {\bf h})= \f{h_{\th}}{r}\p_{\th}(u_{\th} h_{\th})\equiv 0$. Then we obtain the following important maximum principle.
\bl\lab{bernoulli law}
{\it Let $({\bf u}, {\bf h})$ be a axially symmetric smooth solution to (\ref{smhd}) with the form (\ref{asymmetric}). Then we have the following important inequality
\be\lab{bernoulli-inequality}
({\bf u}\cdot\nabla) \Phi- \Delta \Phi\leq 0.
\ee
Hence we have the following maximum principle for $\Phi$ in any bounded domains $\Om$:
\be\lab{max-bounded}
\max_{x\in\Om} \Phi(x)\leq \max_{x\in\p\Om} \Phi(x).
\ee
}
\el
In our cases, as shown in \cite{galdi11}, by adding a constant if necessary, one has $\displaystyle\lim_{|x|\to \oo} p(x)=0$, so $\displaystyle\lim_{|x|\to \oo} \Phi(x)=0$. Hence by maximum principle, we have
\be\lab{max-r3}
\Phi({\bf x})\leq 0,\q \q \forall {\bf x}\in \mbR^3.
\ee

If $({\bf u},{\bf h})$ are axi-symmetric with the form (\ref{asymmetric}), then (\ref{smhd}) can be rewritten as
\be\lab{asym-smhd}
\begin{cases}
(u_r\p_r+ u_z\p_z) u_r-\f{u_{\th}^2}{r} + \p_r p = -\f{h_{\th}^2}{r} + \left(\p_r^2+\f{1}{r}\p_r+\p_z^2-\f{1}{r^2}\right) u_r, \\
(u_r\p_r+ u_z\p_z) u_{\th} +\f{u_r u_{\th}}{r}=\left(\p_r^2+\f{1}{r}\p_r+\p_z^2-\f{1}{r^2}\right)u_{\th},\\
(u_r\p_r+ u_z\p_z) u_z + \p_z p=\left(\p_r^2+\f{1}{r}\p_r+\p_z^2\right) u_z,\\
(u_r\p_r+ u_z\p_z) h_{\th}- \f{u_r h_{\th}}{r}=\left(\p_r^2+\f{1}{r}\p_r+\p_z^2-\f{1}{r^2}\right) h_{\th},\\
\p_r u_r +\f{u_r}{r} + \p_z u_z =0.
\end{cases}
\ee

\bt\lab{liouville-mhd}
{\it Let $({\bf u}, {\bf h})$ be an axially symmetric smooth solution to (\ref{smhd}) with the form (\ref{asymmetric}). 
Then ${\bf h} \equiv 0$.
}
\et
\bpf

Define $\Pi(r,z)=\f{h_{\th}(r,z)}{r}$, then it follows from (\ref{asym-smhd}) that
\be\lab{pi}
(u_r\p_r+ u_z\p_z) \Pi= \left(\p_r^2+\f{1}{r}\p_r+\p_z^2\right) \Pi+\f{2}{r}\p_r\Pi.
\ee
Since $\lim_{|x|\to \oo} \Pi({\bf x})=0$, by the equation (\ref{pi}), we have the maximum and minimum principle, which implies that $\Pi({\bf x})\equiv 0$, i.e. ${\bf h}\equiv 0$. 

\epf

Now we consider the non-resistive, inviscid MHD equations:
\be\label{inviscid-mhd}\begin{cases}
({\bf u}\cdot\nabla) {\bf u}+\nabla p =({\bf h}\cdot\nabla) {\bf h},\q\q\forall x\in \mbR^3,\\
\text{div } {\bf u}=0,\\
({\bf u}\cdot\nabla){\bf h}-({\bf h}\cdot\nabla){\bf u}=0,\q\q\forall x\in \mbR^3,\\
\text{div } {\bf h}=0,\\
\displaystyle\lim_{|x|\to \oo} {\bf u}({\bf x}) =\lim_{|x|\to \oo} {\bf h}({\bf x})= {\bf 0}.
\end{cases}
\ee
For smooth axially symmetric solution $({\bf u}, {\bf h}, p)$ to (\ref{inviscid-mhd}) with the form (\ref{asymmetric}), we can derive the Bernoulli's law for the total head pressure:
\be\lab{bernoulli}
({\bf u}\cdot\na) \Phi=0.
\ee
One can use this Bernoulli's law to establish the existence of weak solutions to steady axially symmetric MHD equations with nonhomogeneous boundary conditions by following the approach developed in \cite{kpr15annals}, for the details please refer to \cite{weng152}.

By slightly modifying the proof in \cite{kpr15}, we can derive the following Liouville theorem for (\ref{inviscid-mhd}).
\bt\lab{liouville-inviscid-mhd}
{\it Let $({\bf u}, {\bf h}, p)$ be an axially symmetric solution to (\ref{inviscid-mhd}) with the form (\ref{asymmetric}) and finite Dirichlet integral. If ${\bf u}$ has no swirl and the corresponding head pressure $\Phi\in L^3(\mbR^3)$ satisfies the condition $\Phi({\bf x})\leq 0$ for any ${\bf x}\in\mbR^3$, then ${\bf u}={\bf h}\equiv 0$.
}
\et

\bpf
Suppose $\na {\bf u},\na {\bf h}\in L^2(\mbR^3)$. Then $({\bf u},{\bf h})\in L^6(\mbR^3)\times L^6(\mbR^3)$. As argued in \cite{kpr15}, we can also derive that
\be\lab{pressure}
p\in D^{2,1}(\mbR^3)\cap D^{1,3/2}(\mbR^3).
\ee
In particular, $\p_r p\in L^1(P_+)$ with $P_+\co \{(r,0,z)\in \mbR^3: r>0, z\in \mbR\}$. Moreover,
\be\lab{pressure-integrability}
p(r,\cdot) \in L^1(\mbR),\q  |{\bf u}(r,\cdot)|^2+ |{\bf h}(r,\cdot)|^2\in L^1(\mbR),\q \q \forall r>0.
\ee

From the steady MHD system (\ref{smhd}) it follows by direct calculation that for any smooth vector function $g$ we have
\be\lab{calculus identity}
\text{div }[p {\bf g}+ ({\bf u}\cdot {\bf g}){\bf u}-({\bf h}\cdot {\bf g}){\bf h}]= p \text{div }{\bf g} + [({\bf u}\cdot \nabla) {\bf g}]\cdot {\bf u}-(({\bf h}\cdot \na{\bf g})\cdot {\bf h}).
\ee

We choose ${\bf g}= g(r) {\bf e}_r$, where ${\bf e}_r$ is the unit vector parallel to the $r$-axis. Then by simple calculations, we obtain
\be\no
&\quad&\text{div }[p {\bf g}+ ({\bf u}\cdot {\bf g}){\bf u}-({\bf h}\cdot {\bf g}){\bf h}]= \text{div }[(p+u_r^2)g(r) {\bf e}_r +g(r) u_r u_z {\bf e}_z] \\\no
&=& \p_r \b(g(r)(p+u_r^2)\b)+\f{g(r)}{r}(p+ u_r^2) +\p_z (g(r) u_r u_z),\\\no
&\quad&  p \text{div }{\bf g} + [({\bf u}\cdot \nabla) {\bf g}]\cdot {\bf u}-(({\bf h}\cdot\na {\bf g})\cdot {\bf h})\\\no
&=& \b(g'(r)+\f{1}{r} g(r)\b)p+ g'(r) u_r^2+\f{g(r)}{r} u_{\th}^2-\f{1}{r} g(r) h_{\th}^2.
\ee

\begin{enumerate}[(i)]
  \item Let $g(r)= r $. Then for axially symmetric ${\bf u}$ and $p$ we get
  \be\no
  p \text{div }{\bf g} + [({\bf u}\cdot \nabla) {\bf g}]\cdot {\bf u}-(({\bf h}\cdot\na {\bf g})\cdot {\bf h})= 2p + u_{\theta}^2 + u_r^2- h_{\th}^2.
  \ee
  Integrating this identity over the three dimensional infinite cylinder $C_t=\{(x_1,x_2,x_3)\in \mbR^3: r=\sqrt{x_1^2+x_2^2}<t, x_3\in \mbR\}$, we obtain
  \be\lab{integral id1}
  t^2 \int_{\mbR} [p(t,z)+ u_r^2(t,z)]dz &=& \iint_{P_t} r[2p +u_{\th}^2+ u_r^2-h_{\th}^2] dr dz\\\no
  &=& \iint_{P_t} r(2\Phi- u_z^2-h_{\th}^2) dr dz \leq 0
  \ee
  where $P_t=\{(r,z)\in P_+: r<t\}$. Here we use $\Phi(x)\leq 0$ for all $x\in \mbR^3$.
  \item Let $g(r)=\f{1}{r}$. Then
  \be\no
  p \text{div }{\bf g} + [({\bf u}\cdot \nabla) {\bf g}]\cdot {\bf u}-(({\bf h}\cdot\na {\bf g})\cdot {\bf h})=\f{1}{r^2} (u_{\th}^2-u_r^2-h_{\th}^2).
  \ee
\end{enumerate}

Since we have an essential singularity at $r=0$, we need to integrate this identity over the cylindrical annulus $C_{t_0t} =\{(x_1,x_2,x_3)\in \mbR^3: r=\sqrt{x_1^2+x_2^2}\in (t_0,t), x_3\in \mbR\}$ to obtain
\be\lab{integral id2}
\int_{\mbR} [p(t,z)+ u_r^2(t,z)]dz- \int_{\mbR} [p(t_0,z)+ u_r^2(t_0,z)] dz =\iint_{P_{t_0 t}} \f1r(u_{\th}^2-u_r^2-h_{\th}^2) dz dr,
\ee
where $P_{t_0t}= \{(r,z)\in P_+: r\in (t_0,t), z\in \mbR\}$.

Since $\int_{\mbR} [p(t,z)+ u_r^2(t,z)] dz\to 0$ as $t\to +\oo$ and
\be\no
\iint_{P_+} \b|\f{1}{r} (u_{\th}^2- u_r^2-h_{\th}^2)\b| dz dr<\oo,
\ee
we obtain immediately from the above formulas that
\be\lab{integral id3}
\int_{\mbR} [p(t,z)+ u_r^2(t,z)] dz &=&\iint_{P_{t\oo}} \left(\f{1}{r}(h_{\th}^2+u_r^2-u_{\theta}^2)\right) dz dr\\\no
&\geq& 0,\q\q \text{if } u_{\th}(r,z)\equiv 0,
\ee
where $P_{t\oo}=\{(r,z)\in P_+: r\in (t,+\oo), z\in \mbR\}$.

If ${\bf u}$ is axially symmetric with no swirl, then the formulas (\ref{integral id1}) and (\ref{integral id3}) imply ${\bf u}\equiv 0$.
Indeed, from the last two inequalities it follows that $\int_{\mbR} [p(t,z)+ u_r^2(t,z)] dz\equiv 0$ for all $t>0$, and thus, by virtue of (\ref{integral id3}) and $u_{\theta}(x)=0$, we obtain $u_r\equiv h_{\th}\equiv 0$. Therefore, by (\ref{inviscid-mhd}) we conclude that $u_z\equiv 0$.

\epf

\section{Liouville theorem for steady viscous resistive Hall-MHD equations}\lab{lhmhd}\hspace*{\parindent}
The steady Hall-MHD equations read as follows.
\be\label{shmhd}\begin{cases}
({\bf u}\cdot\nabla) {\bf u}+\nabla p =({\bf h}\cdot\nabla) {\bf h}+ \Delta {\bf u},\q\q\forall x\in \mbR^3,\\
\text{div } {\bf u}=0,\\
({\bf u}\cdot\nabla){\bf h}-({\bf h}\cdot\nabla){\bf u}+\nabla \times ((\nabla\times {\bf h})\times {\bf h})=\Delta {\bf h},\q\q\forall x\in \mbR^3,\\
\text{div } {\bf h}=0,\\
\displaystyle\lim_{|x|\to \oo} {\bf u}({\bf x}) =\lim_{|x|\to \oo} {\bf h}({\bf x})= {\bf 0}.
\end{cases}
\ee

In this section, we consider the smooth solution to (\ref{shmhd}) with finite Dirichlet integral (\ref{mhd-dirichlet}). Comparing with the well-known MHD system, the Hall term $\nabla\times ((\nabla\times {\bf h})\times {\bf h})$ is included due to the Ohm's law, which is believed to be a key issue for understanding magnetic reconnection. Note that the Hall term is quadratic in the magnetic field and involves the second order derivatives. A derivation of Hall-MHD system from a two-fluids Euler-Maxwell system for electrons and ions was presented in \cite{adfl}, through a set of scaling limits. They also provided a kinetic formulation for the Hall-MHD, and proved the existence of global weak solutions for the incompressible viscous resistive Hall-MHD system. The authors in \cite{cweng} have showed unsteady Hall-MHD system without resistivity may develop finite time singularity for a special class of axially symmetric datum. Different from the steady MHD case, a weak solution to (\ref{shmhd}) with finite Dirichlet integral may not be smooth. Chae and Wolf \cite{cwolf15} have investigated the partial regularity of suitable weak solutions to steady Hall-MHD equations, showing the set of possible singularities has Hausdorff dimension at most one. The authors in \cite{cdl} had derived a Liouville theorem for steady Hall-MHD under the conditions $({\bf u}, {\bf h})\in L^{\oo}(\mbR^3)\cap L^{\f 92}(\mbR^3)$.

\bt\lab{Liou-hmhd}
{\it Let $({\bf u}, {\bf h})$ be a smooth solution to (\ref{shmhd}) in $\mbR^3$ with finite Dirichlet integral. If ${\bf u}\in L^3(\mbR^3)$, then ${\bf u}= {\bf h}\equiv 0$.
}
\et
\bpf
Similar to the proof in Theorem \ref{Liou-mhd}, we only need to check the effect of the Hall-term $H=\int_{\mbR^3} \si_R(x) {\bf h}\cdot (\nabla\times ((\nabla \times{\bf h})\times {\bf h})) dx$. Indeed,
\be\no
|H|&=& \b|-\int_{\mbR^3} \f{1}{R} {\bf h}\cdot \left(\nabla \si(x/R)\times ((\nabla\times {\bf h})\times {\bf h})\right) dx\b| \\\no
&\leq& \|\nabla {\bf h}\|_{L^2(Q_R)}\|{\bf h}\|_{L^6(Q_R)}^2 \to 0,\q \text{as }R\to \oo.
\ee
Then we finish the proof.
\epf

If we also consider the axially symmetric smooth solution $({\bf u}, {\bf h})$ to (\ref{shmhd}) with the form
\be\lab{special-hmhd}
{\bf u}({\bf x})=u_r(r,z) {\bf e}_r+ u_{\th}(r,z) {\bf e}_{\th}+ u_z(r,z) {\bf e}_z,\q {\bf h}({\bf x})= h_{\th}(r,z) {\bf e}_{\th},
\ee
then (\ref{shmhd}) can be rewritten as
\be\lab{asym-shmhd}
\begin{cases}
(u_r\p_r+ u_z\p_z) u_r-\f{u_{\th}^2}{r} + \p_r p = -\f{h_{\th}^2}{r} + \left(\p_r^2+\f{1}{r}\p_r+\p_z^2-\f{1}{r^2}\right) u_r, \\
(u_r\p_r+ u_z\p_z) u_{\th} +\f{u_r u_{\th}}{r}=\left(\p_r^2+\f{1}{r}\p_r+\p_z^2-\f{1}{r^2}\right)u_{\th},\\
(u_r\p_r+ u_z\p_z) u_z + \p_z p=\left(\p_r^2+\f{1}{r}\p_r+\p_z^2\right) u_z,\\
(u_r\p_r+ u_z\p_z) h_{\th}- \f{u_r h_{\th}}{r}- \f{2 h_{\th}}{r}\p_z h_{\th}=\left(\p_r^2+\f{1}{r}\p_r+\p_z^2-\f{1}{r^2}\right) h_{\th},\\
\p_r u_r +\f{u_r}{r} + \p_z u_z =0.
\end{cases}
\ee
Due to Hall term, we can not arrive at any maximum principle for the total head pressure. However, we still have the following Liouville theorem.
\bt\lab{Liouville-HMHD}
{\it Let ${\bf u}$ and ${\bf h}$ be an axially symmetric smooth solution to (\ref{shmhd}) with the form (\ref{special-hmhd}). 
then ${\bf h}\equiv 0$.
}
\et

\bpf

Setting $\Pi\co \f{h_{\th}}{r}$ as before, then $\Pi$ satisfies the following equation
\be\lab{pi-hmhd}
(u_r\p_r+u_z\p_z)\Pi- 2\Pi\p_z \Pi= \left(\p_r^2+\f{1}{r}\p_r +\p_z^2\right)\Pi + \f{2}{r}\p_r \Pi.
\ee

Since $\displaystyle \lim_{|x|\to \oo} {\bf h}({\bf x})=0$, then $\displaystyle\lim_{|x|\to \oo} \Pi({\bf x})=0$. Then by the maximum principle of $\Pi$ from (\ref{pi-hmhd}), we conclude that $\Pi\equiv 0$, i.e. $h_{\th}\equiv 0$. 

\epf

{\bf Acknowledgement.}  Weng's research was supported by Basic Science Research Program through the National Research Foundation of Korea(NRF) funded by the Ministry of Education, Science and Technology (2015049582).


\end{document}